\documentclass[secnum,leqno]{rims-bessatsu}%%% eqno resets every section 
\usepackage{amssymb, amsmath}%%% option packages
\usepackage{graphics}%%% option packages
\usepackage[dvips]{graphicx}%%% option packages
\usepackage{eepic}%%% option packages
\usepackage{epic}%%% option packages
\usepackage{color}

\definecolor{blue}{rgb}{0,0,1}
\definecolor{red}{rgb}{1,0,0}
\definecolor{green}{rgb}{0,1,0}
\long\def\red#1\endred{\textcolor{red}{#1}}
\long\def\blue#1\endblue{\textcolor{blue}{#1}}
\long\def\purple#1\endpurple{\textcolor{purple}{
{\raisebox{1.6ex}{\makebox[0cm][l]{\tiny ch}}}#1
{\raisebox{1.6ex}{\makebox[0cm][r]{\tiny endch}}}}}
\long\def\green#1\endgreen{\textcolor{green}{#1}}
\newtheorem{thm}{Theorem}[section]

\theoremstyle{definition}

\theoremstyle{remark}
\newtheorem{rem}{Remark}%%%
\title{Functional equations for double series of Euler-Hurwitz-Barnes type
with coefficients}
\TitleHead{Functional equations for double series}          %optional
\author{\textsc{YoungJu Choie}\footnote{Department of Mathematics and PMI,
Pohang University of Science and Technology, Pohang 790-784, Korea.\newline e-mail: 
\texttt{yjc@postech.ac.kr}}
          ~and \textsc{Kohji Matsumoto}\footnote{Graduate School of Mathematics, 
Nagoya University, Chikusa-ku, Nagoya
464-8602, Japan. \endgraf e-mail: \texttt{kohjimat@math.nagoya-u.ac.jp}}}
\AuthorHead{Choie and Matsumoto}         %optional but recommended
\support{Supported by ***}  %optional
         %optional
\VolumeNo{x}           %editors must define this.
\YearNo{201x}           %editors must define this.
\PagesNo{000--000}      %editors must define this.
\communication{Received April 20, 201x.
Revised September 11, 201x.}      %editors must define this.
\begin{document}
%
% The text goes here.  
% Be sure to use the appropriate "theorem-like" environment as 
% is the following examples.  Never use plain TeX commands for these, as
% they will cause interference with the styles of other papers. 

\maketitle

%\tableofcontents      %optional
\begin{abstract}      %optional
We first survey the known results on functional equations for the double
zeta-function of Euler type and its various generalizations.   Then we prove
two new functional equations for double series of Euler-Hurwitz-Barnes type
with complex coefficients.   The first one is of general nature, while the second
one is valid when the coefficients are Fourier coefficients of a cusp form.
\end{abstract}
%%%%%%%%%%%%%%%%%%%%%%%%%%%%%%%%%%%%%%%%%%%%%%%%%%%%%%%%%%%%%%%%%%%%%%%%%%%%%%%
\section{Introduction}
%%%%%%%%%%%%%%%%%%%%%%%%%%%%%%%%%%%%%%%%%%%%%%%%%%%%%%%%%%%%%%%%%%%%%%%%%%%%%%%

Functional equations play a very important role in the theory of zeta and
$L$-functions.   In the case of the most classical Riemann zeta-function
$\zeta(s)$, Riemann proved the beautiful symmetric functional equation
\begin{align}\label{1-1}
\pi^{-s/2}\Gamma\left(\frac{s}{2}\right)\zeta(s)=
\pi^{-(1-s)/2}\Gamma\left(\frac{1-s}{2}\right)\zeta(1-s).
\end{align}
This symmetry, however, does not remain when we consider the Hurwitz
zeta-function $\zeta(s,\alpha)=\sum_{n=0}^{\infty}(n+\alpha)^{-s}$, where
$0<\alpha\leq 1$.   The functional equation for $\zeta(s,\alpha)$ is of the form
\begin{align}\label{1-2}
\zeta(s,\alpha)=\frac{\Gamma(1-s)}{i(2\pi)^{1-s}}\left\{
e^{\pi is/2}\phi(1-s,\alpha)-e^{-\pi is/2}\phi(1-s,-\alpha)\right\},
\end{align}
where $i=\sqrt{-1}$ and $\phi(s,\alpha)=\sum_{n=1}^{\infty}e^{2\pi in\alpha}
n^{-s}$ is the Lerch zeta-function (see \cite[(2.17.3)]{Tit51}).
When $\alpha=1$, \eqref{1-2} reduces to \eqref{1-1}, but the form \eqref{1-2} is
no longer symmetric.

In recent decades, the theory of various multiple zeta-functions has been 
developed extensively.   Therefore it is natural to search for 
functional equations for those multiple zeta-functions.   For example, Barnes
multiple zeta-functions
\begin{align}\label{1-3}
\sum_{n_1=1}^{\infty}\cdots\sum_{n_r=1}^{\infty}(\omega_1 n_1+\cdots+
\omega_r n_r+\alpha)^{-s}
\end{align}
(Barnes \cite{Bar04}) with complex parameters $\omega_1,\ldots,\omega_r,\alpha$
is a direct generalization of Hurwitz zeta-functions, so a kind of functional
equation similar to \eqref{1-2} is expected to hold.   In fact, such an equation 
has been proved under certain condition of parameters (see Hardy and Littlewood
\cite{HL22} in the case $r=2$, and \cite{KMT-JRMS} \cite{Shi13} in general).

Another important class of multiple zeta-functions is the multi-variable sum
\begin{align}\label{1-4}
\zeta_r(s_1,\ldots,s_r)=\sum_{n_1=1}^{\infty}\cdots\sum_{n_r=1}^{\infty}
n_1^{-s_1}(n_1+n_2)^{-s_2}\cdots(n_1+\cdots+n_r)^{-s_r},
\end{align}
which was first considered by Euler in the case $r=2$, and then introduced by
Hoffman \cite{Hof92} and Zagier \cite{Zag94} independently of each other for
general $r$.
To find some kind of functional equations for the sum \eqref{1-4} or its
variants/generalizations seems a rather complicated problem.   Let us consider
the simplest case $r=2$:
\begin{align}\label{1-5}
\zeta_2(s_1,s_2)=\sum_{m=1}^{\infty}\sum_{n=1}^{\infty}m^{-s_1}(m+n)^{-s_2}.
\end{align}
In the following sections we will report our attempt to search for functional
equations for \eqref{1-5} and its various variants.

%%%%%%%%%%%%%%%%%%%%%%%%%%%%%%%%%%%%%%%%%%%%%%%%%%%%%%%%%%%%%%%%%%%%%%%%%%%%%%%%
\section{Functional equations for double zeta-functions}
%%%%%%%%%%%%%%%%%%%%%%%%%%%%%%%%%%%%%%%%%%%%%%%%%%%%%%%%%%%%%%%%%%%%%%%%%%%%%%%%

Let $0<\alpha\leq 1$, $0\leq \beta\leq 1$, $\omega>0$, and define 
\begin{align}\label{2-1}
\zeta_2(s_1,s_2;\alpha,\beta,\omega)=\sum_{m=0}^{\infty}\sum_{n=1}^{\infty}
\frac{e^{2\pi in\beta}}{(\alpha+m)^{s_1}(\alpha+m+n\omega)^{s_2}}.
\end{align}
This is a generalization of \eqref{1-5} of the Hurwitz-Lerch type.
The functional equation for $\zeta_2(s_1,s_2;\alpha,\beta,\omega)$ can be stated
in terms of the confluent hypergeometric function
\begin{align}\label{2-1a}
\Psi(a,c;x)=\frac{1}{\Gamma(a)}\int_0^{e^{i\phi}\infty}e^{-xy}y^{a-1}
(1+y)^{c-a-1}dy,
\end{align}
where $a,c,x\in\mathbb{C}$, $\Re a>0$, $-\pi<\phi<\pi$, $|\phi+\arg x|<\pi/2$.
Define
\begin{align}\label{2-2}
F_{\pm}(s_1,s_2;\alpha,\beta,\omega)=\sum_{k=1}^{\infty}
\sigma_{s_1+s_2-1}(k;\alpha,\beta)\Psi(s_2,s_1+s_2;\pm 2\pi ik\omega),
\end{align}
where
$$
\sigma_c(k;\alpha,\beta)=\sum_{0<d|k}e^{2\pi id\alpha}e^{2\pi i(k/d)\beta}d^c.
$$
Then the following functional equation is known.

\begin{thm}\label{thm1}
The functions $F_{\pm}(s_1,s_2;\alpha,\beta,\omega)$ can be continued
meromorphically to the whole space $\mathbb{C}^2$, and the functional equation
\begin{align}\label{2-3}
\lefteqn{\zeta_2(s_1,s_2;\alpha,\beta,\omega)
=\frac{\Gamma(1-s_1)}{\Gamma(s_2)}\Gamma(s_1+s_2-1)\phi(s_1+s_2-1,\beta)
  \omega^{1-s_1-s_2}}\\
&\;+\Gamma(1-s_1)\omega^{1-s_1-s_2}\notag\\
&\quad\times\{F_+(1-s_2,1-s_1;\beta,\alpha,\omega)+
F_-(1-s_2,1-s_1;\beta,-\alpha,\omega)\}\notag
\end{align}
holds.
\end{thm}

Here, the first term on the right-hand side of \eqref{2-3} is an ``additional'' 
term, and the main body of the right-hand side is the second term involving
$F_{\pm}$.   This term, compared with the left-hand side,
expresses a duality between the values at $(s_1,s_2)$ and
at $(1-s_2,1-s_1)$ (and also a duality between $(\alpha,\beta)$ and
$(\beta,\alpha)$), and hence formula \eqref{2-3} can be regarded as a double 
analogue of \eqref{1-2}.   The functions $F_{\pm}$ are not Dirichlet series in
the usual sense, but the confluent hypergeometric function satisfies the
asymptotic expansion 
\begin{align}\label{2-4}
\Psi(a,c;x)=x^{-a}-a(a-c+1)x^{-a-1}+\frac{a(a+1)(a-c+1)(a-c+2)}{2}x^{-a-2}
+\cdots,
\end{align}
so we may say that $\Psi(a,c;x)$ can be approximated by $x^{-a}$.   From this
viewpoint it is possible to say that $F_{\pm}$ can be approximated by the
Dirichlet series
$$
\sum_{k=1}^{\infty}\sigma_{s_1+s_2-1}(k;\alpha,\beta)(\pm 2\pi ik\omega)^{-s_2}.
$$
In this sense $F_{\pm}$ may be regarded as generalized Dirichlet series.

The meromorphic continuation of $F_{\pm}$ was shown in \cite[Proposition 2]{Mat04},
where their functional equation
\begin{align}\label{2-5}
F_{\pm}(1-s_2,1-s_1;\beta,\alpha,\omega)=
(\pm 2\pi i\omega)^{s_1+s_2-1}F_{\pm}(s_1,s_2;\alpha,\beta,\omega) 
\end{align}
was also proved.   The transformation formula
\begin{align}\label{2-6}
\Psi(a,c;x)=x^{1-c}\Psi(a-c+1,2-c;x)
\end{align}
of the confluent hypergeometric function is used in the proof of \eqref{2-5}.
Applying \eqref{2-5} to \cite[Proposition 1]{Mat04}, we can immediately obtain
formula \eqref{2-3}.   Therefore the above Theorem \ref{thm1} is essentially
included in \cite{Mat04}, though it is first explicitly stated in \cite{Mat07}.

The main statement of \cite{Mat04} is as follows.   Let
\begin{align}\label{2-7}
\lefteqn{g(s_1,s_2;\alpha,\beta,\omega)=
\zeta_2(s_1,s_2;\alpha,\beta,\omega)}\\                                             
&-\frac{\Gamma(1-s_1)}{\Gamma(s_2)}\Gamma(s_1+s_2-1)\phi(s_1+s_2-1,\beta)
  \omega^{1-s_1-s_2}\notag.
\end{align}
Then

\begin{thm}[{\cite[Theorem 2]{Mat04}}]
\label{thm2}
We have
\begin{align}\label{2-8}
\lefteqn{\frac{g(s_1,s_2;\alpha,\beta,\omega)}{(2\pi)^{s_1+s_2-1}\Gamma(1-s_1)}
=\frac{g(1-s_2,1-s_1;1-\beta,1-\alpha,\omega)}{(i\omega)^{s_1+s_2-1}\Gamma(s_2)}}\\
&+e^{\pi i(s_1+s_2-1)/2}F_+(s_1,s_2;\alpha,\beta,\omega)
-e^{\pi i(1-s_1-s_2)/2}F_+(s_1,s_2;1-\alpha,1-\beta,\omega),\notag
\end{align}
and especially, when $\alpha=\beta=1$, we have
\begin{align}\label{2-9}
\lefteqn{\frac{g(s_1,s_2;1,1,\omega)}{(2\pi)^{s_1+s_2-1}\Gamma(1-s_1)}        
=\frac{g(1-s_2,1-s_1;1,1,\omega)}{(i\omega)^{s_1+s_2-1}\Gamma(s_2)}}\\
&+2i\sin\left(\frac{\pi}{2}(s_1+s_2-1)\right)F_+(s_1,s_2;1,1,\omega).\notag
\end{align}
\end{thm}

This theorem can also be easily deduced from Proposition 1 and Proposition 2 of 
\cite{Mat04}.

\begin{rem}\label{rem1} (Historical note)
The idea of the proof of Theorem \ref{thm1} goes back to \cite{KM91}, where  
certain mean values of Dirichlet $L$-functions were studied.
The application of the confluent hypergeometric function in this context was
first done by Katsurada \cite{Kat93}.   In order to study Barnes' double 
zeta-functions (the case $r=2$ of \eqref{1-3}), the second-named author \cite{Mat98}
introduced the two-variable double series
\begin{align}\label{2-10}
\zeta_2(s_1,s_2;\alpha,\omega)=\sum_{m=0}^{\infty}\sum_{n=1}^{\infty}               
{(\alpha+m)^{ {-}  s_1}(  \alpha+m+n\omega)^{ {-}   s_2}}
\end{align}
(the case $\beta=0$ of \eqref{2-1}),
and studied its properties, invoking the methods in \cite{KM91}, \cite{Kat93}.
In particular, the special case $\beta=0$ of \cite[Proposition 1]{Mat04} was
already given in \cite[(5.5)]{Mat98}.
\end{rem}

%%%%%%%%%%%%%%%%%%%%%%%%%%%%%%%%%%%%%%%%%%%%%%%%%%%%%%%%%%%%%%%%%%%%%%%%%%%%%%%%
\section{The symmetric form}
%%%%%%%%%%%%%%%%%%%%%%%%%%%%%%%%%%%%%%%%%%%%%%%%%%%%%%%%%%%%%%%%%%%%%%%%%%%%%%%%

In this and the next section 
we present the contents of two joint papers of Komori, Tsumura
and the second-named author.   Theorem \ref{thm1} is a non-symmetric functional
equation, similarly to \eqref{1-2} for Hurwitz zeta-functions.   Is there any 
{\it symmetric} functional equation for double zeta-functions?   One of the main
point of \cite{KMT-Debrecen} is that such equations do exist on certain special
hyperplanes.

In \cite{KMT-Debrecen}, the following generalization of \eqref{1-5} was introduced:
\begin{align}\label{3-1}
\zeta_2(s_1,s_2;\omega_1,\omega_2)=\sum_{m=1}^{\infty}\sum_{n=1}^{\infty}
(m\omega_1)^{-s_1}(m\omega_1+n\omega_2)^{-s_2},
\end{align}
where $\omega_1,\omega_2\in\mathbb{C}$ with $\Re\omega_1>0$, $\Re\omega_2>0$.
The functional equation similar to Theorem \ref{thm1} holds for 
$\zeta_2(s_1,s_2;\omega_1,\omega_2)$, which is
\begin{align}\label{3-2}
\lefteqn{\zeta_2(s_1,s_2;\omega_1,\omega_2)                                           
=\frac{\Gamma(1-s_1)}{\Gamma(s_2)}\Gamma(s_1+s_2-1)\zeta(s_1+s_2-1)\omega_1^{-1}
  \omega_2^{1-s_1-s_2}}\\ 
&\;+\Gamma(1-s_1)\omega_1^{-1}\omega_2^{1-s_1-s_2}\notag\\  
&\quad\times\{F_+(1-s_2,1-s_1;1,1,\omega_2/\omega_1)+
F_-(1-s_2,1-s_1;1,1,\omega_2/\omega_1)\}\notag                                           
\end{align}
(\cite[Theorem 2.1]{KMT-Debrecen}).   This formula itself can be proved just
similarly to Theorem \ref{thm1}.   However, from this formula, it is possible to 
deduce the following symmetric functional equation.   Let
\begin{align}\label{3-3}
\lefteqn{\xi(s_1,s_2;\omega_1,\omega_2)=\left(\frac{2\pi i}{\omega_1\omega_2}\right)
^{(1-s_1-s_2)/2}\Gamma(s_2)}\\
&\times\left\{\zeta_2(s_1,s_2;\omega_1,\omega_2)-\frac{\Gamma(1-s_1)}{\Gamma(s_2)}
\Gamma(s_1+s_2-1)\zeta(s_1+s_2-1)\omega_1^{-1}\omega_2^{1-s_1-s_2}\right\}.\notag
\end{align}
Then we have

\begin{thm}[{\cite[Theorem 2.2]{KMT-Debrecen}}]
\label{thm3}
The hyperplane
\begin{align}\label{3-4}
\Omega_{2k+1}=\{(s_1,s_2)\in\mathbb{C}^2\;|\;s_1+s_2=2k+1\}
\end{align}
$(k\in\mathbb{Z}\setminus\{0\})$ is not a singular locus of
$\zeta_2(s_1,s_2;\omega_1,\omega_2)$, and when $(s_1,s_2)\in \Omega_{2k+1}$, the 
functional equation
\begin{align}\label{3-5}
\xi(s_1,s_2;\omega_1,\omega_2)=\xi(1-s_2,1-s_1;\omega_1,\omega_2)
\end{align}
holds.
\end{thm}

When $\omega_1=1$ and $\omega_2=\omega$, this theorem is actually almost
immediately obtained from \eqref{2-9}.   (But the second-named author did not notice
this point when he wrote \cite{Mat04}).

From Theorem \ref{thm3}, we can evaluate certain values of
$\zeta_2(s_1,s_2;\omega_1,\omega_2)$ when $s_1,s_2$ are non-positive integers
(\cite[Corollary 2.4]{KMT-Debrecen}).

%%%%%%%%%%%%%%%%%%%%%%%%%%%%%%%%%%%%%%%%%%%%%%%%%%%%%%%%%%%%%%%%%%%%%%%%%%%%%%%
\section{Functional equations for double $L$-functions}
%%%%%%%%%%%%%%%%%%%%%%%%%%%%%%%%%%%%%%%%%%%%%%%%%%%%%%%%%%%%%%%%%%%%%%%%%%%%%%%

In the statement of aforementioned functional equations (\eqref{2-3}, 
\eqref{3-2}), an ``additional'' term involving
$\phi(s_1+s_2-1,\beta)$ or $\zeta(s_1+s_2-1)$ appears as the first term of the
right-hand side, which also appears in the definitions \eqref{2-7} and \eqref{3-3}.
The reason of the appearance of such a term can be clarified when we consider a 
little more general situation.

Let $f\in\mathbb{Z}_{\geq 2}$, and let $a_j:\mathbb{Z}\to\mathbb{C}$ ($j=1,2$)
be periodic functions with period $f$.   Assume $a_j$ is an even (or odd) function,
and let $\lambda(a_j)=1$ (resp.$\;-1$) if $a_j$ is even (resp. odd).   Define the
associated $L$-function by $L(s,a_j)=\sum_{m=1}^{\infty}a_j(m)m^{-s}$.
In \cite{KMT-IJNT} the double series
\begin{align}\label{4-1}
L_2(s_1,s_2;a_1,a_2;\omega_1,\omega_2)=\sum_{m=1}^{\infty}\sum_{n=1}^{\infty}
\frac{a_1(m)a_2(n)}{(m\omega_1)^{s_1}(m\omega_1+n\omega_2)^{s_2}}
\end{align}
was introduced.   Denote the finite Fourier expansion of $a_j$ by
$$
a_j(m)=\sum_{\nu=1}^f \widehat{a}_j(\nu)e^{2\pi i\nu m/f}.
$$
The series $L_2(s_1,s_2;a_1,a_2;\omega_1,\omega_2)$ can be continued meromorphically 
to the whole space $\mathbb{C}^2$.   As an analogue of Theorem \ref{thm3}, the 
following symmetric functional equation holds.

\begin{thm}[{\cite[Theorem 2.1]{KMT-IJNT}}]
\label{thm4}
The functional equation
\begin{align}\label{4-2}
&\left(\frac{2\pi i}{f\omega_1\omega_2}\right)^{(1-s_1-s_2)/2}
\Biggl\{\Gamma(s_2)L_2(s_1,s_2;a_1,a_2;\omega_1,\omega_2)\Biggr.\\
&\quad\Biggl.-\frac{\omega_2^{1-s_1-s_2}}{f\omega_1}\Gamma(1-s_1)\Gamma(s_1+s_2-1)
L(s_1+s_2-1,a_2)\sum_{\nu=1}^f a_1(\nu)\Biggr\}\notag\\
&=\left(\frac{2\pi i}{f\omega_1\omega_2}\right)^{(s_1+s_2-1)/2}
\Biggl\{\Gamma(1-s_1)L_2(1-s_2,1-s_1;\widehat{a}_2,\widehat{a}_1;\omega_1,\omega_2)
\Biggr.\notag\\
&\quad\Biggl.-\frac{\omega_2^{s_1+s_2-1}}{f\omega_1}\Gamma(s_2)\Gamma(1-s_1-s_2)
L(1-s_1-s_2,\widehat{a}_1)\sum_{\nu=1}^f \widehat{a}_2(\nu)\Biggr\}\notag
\end{align}
holds on the hyperplane $s_1+s_2=2k+1$ $(k\in\mathbb{Z})$ if 
$\lambda(a_1)\lambda(a_2)=1$, and on the hyperplane $s_1+s_2=2k$ $(k\in\mathbb{Z})$ 
if $\lambda(a_1)\lambda(a_2)=-1$. 
\end{thm} 

Therefore if
\begin{align}\label{4-3}
\sum_{\nu=1}^f a_1(\nu)=\sum_{\nu=1}^f \widehat{a}_2(\nu)=0
\end{align}
then the ``additional'' terms do not appear.
In particular, if $a_1, a_2$ are primitive Dirichlet characters $\chi_1,\chi_2$
of conductor $f$, then \eqref{4-3} holds, and so the functional equation becomes
the following very simple form:
\begin{align}\label{4-4}                                                                
&\left(\frac{2\pi i}{f\omega_1\omega_2}\right)^{(1-s_1-s_2)/2}                          
\Gamma(s_2)L_2(s_1,s_2;\chi_1,\chi_2;\omega_1,\omega_2)\\
&=\left(\frac{2\pi i}{f\omega_1\omega_2}\right)^{(s_1+s_2-1)/2}                         
\Gamma(1-s_1)L_2(1-s_2,1-s_1;\widehat{\chi}_2,\widehat{\chi}_1;\omega_1,\omega_2)      
\notag                                                                         
\end{align}
on the hyperplanes indicated in the statement of Theorem \ref{thm4}.

\begin{rem}\label{rem2}
Using the relation 
$\widehat{\chi}_j(m)=\overline{\chi}_j(m)/\tau(\overline{\chi}_j)$ (where
$\tau(\cdot)$ denotes the Gauss sum), we can restate formula \eqref{4-4} in which 
$\widehat{\chi}_j$ is replaced by $\overline{\chi}_j$.   Such a formula is
stated as \cite[Corollary 2.3]{KMT-IJNT}, but the factor $\chi_1(-1)f$ is lacking
on the left-hand side of the statement there.
\end{rem}

From this result we can evaluate certain values of double $L$-functions at
non-positive integer points (\cite[Section 3]{KMT-IJNT}).
Those results further motivates the construction of the theory of $p$-adic
multiple $L$-functions.   A double analogue of the Kubota-Leopoldt $p$-adic
$L$-function was introduced in \cite[Section 4]{KMT-IJNT}, and then, a more
general theory of $p$-adic multiple $L$-functions has been developed in
\cite{FKMT}.

To prove Theorem \ref{thm4}, the following double zeta-function of Hurwitz-Lerch
type was introduced:
\begin{align}\label{4-5}
\zeta_2(s_1,s_2;\alpha,\beta;\omega_1,\omega_2)=\sum_{m=0}^{\infty}
\sum_{n=1}^{\infty}\frac{e^{2\pi in\beta}}
{((\alpha+m)\omega_1)^{s_1}((\alpha+m)\omega_1+n\omega_2)^{s_2}},
\end{align}
where $0<\alpha\leq 1$, $0\leq\beta\leq 1$, $\omega_1,\omega_2\in\mathbb{C}$ with
$\Re\omega_1>0$, $\Re\omega_2>0$.   This is a generalization of both \eqref{2-1} 
and \eqref{3-1}.   
(Frankly speaking, this is almost equivalent to consider \eqref{2-1} with
$\omega\in\mathbb{C}$, $|\arg\omega|<\pi$.)
Since $L_2(s_1,s_2;a_1,a_2;\omega_1,\omega_2)$ can be written
as a linear combination of functions of the form \eqref{4-5}, it is sufficient to
show a functional equation for \eqref{4-5}.   This can be done in a way similar
to the proofs of Theorem \ref{thm1} and Theorem \ref{thm3} (see 
\cite[Sections 5, 6]{KMT-IJNT}).

%%%%%%%%%%%%%%%%%%%%%%%%%%%%%%%%%%%%%%%%%%%%%%%%%%%%%%%%%%%%%%%%%%%%%%%%%%%%%%%%
\section{Functional equations for double series with complex coefficients}
%%%%%%%%%%%%%%%%%%%%%%%%%%%%%%%%%%%%%%%%%%%%%%%%%%%%%%%%%%%%%%%%%%%%%%%%%%%%%%%%

Let ${\mathfrak A}=\{a(n)\}_{n\geq 1}$ be a sequence of complex numbers.
In \cite{CM}, the authors considered the double series
\begin{align}\label{5-1}
L_2(s_1,s_2;{\mathfrak A})=\sum_{m=1}^{\infty}\sum_{n=1}^{\infty}
\frac{a(n)}{m^{s_1}(m+n)^{s_2}}.
\end{align}
The original motivation of \cite{CM} is to study the case when $a(n)$'s are
Fourier coefficients of modular forms.
In \cite{Man05} \cite{Man06}, Manin extended the theory of periods of modular forms
replacing integration along geodesics by iterated integrations.   In this context
some kind of multiple Dirichlet series with Fourier coefficients on the numerator
naturally appears (cf. \cite{CI}).   
Multiple series with Fourier coefficients is also expected 
to be useful to evaluate certain multiple sums of Fourier coefficients
(analogously to de la Bret{\`e}che \cite{Bre01}, Ishikawa \cite{Ish02},
Essouabri \cite{Ess12} etc).

In the present paper we consider a slight generalization of \eqref{5-1}, that is
\begin{align}\label{5-2}
L_2(s_1,s_2;\alpha;\omega;{\mathfrak A})=
\sum_{m=0}^{\infty}\sum_{n=1}^{\infty}\frac{a(n)}
{(\alpha+m)^{s_1}(\alpha+m+n\omega)^{s_2}},
\end{align}
where $\omega\in\mathbb{C}$ with $|\arg\omega|<\pi$.
This is a generalization of both \eqref{2-1} and \eqref{5-1}.  (The exponential factor
$e^{2\pi in\beta}$ on the numerator of \eqref{2-1} is to be included in $a(n)$.)
The following two theorems (Theorems \ref{thm5} and \ref{thm6}) are the main results
of the present paper.   When $\alpha=\omega=1$, those theorems are
Theorem 1.1 and Theorem 1.5 in \cite{CM}, respectively.

Concerning the order of the coefficients $a(n)$, we assume:

(i) $a(n)\ll n^{(\kappa-1)/2+\varepsilon}$ with a certain constant $\kappa>1$,
where $\varepsilon$ is an arbitrarily small positive number.

Then, the Dirichlet series
\begin{align}\label{5-3}
L(s,{\mathfrak A})= \sum_{n=1}^{\infty}a(n)n^{-s}
\end{align}
is absolutely convergent in the region $\Re s>(\kappa+1)/2$, and we further assume

(ii) \eqref{5-3} can be continued meromorphically to the whole complex plane with
only finitely many poles.

Under these two assumptions, we see that \eqref{5-2} is absolutely convergent
in the region
\begin{align}\label{5-4}
\Re s_2>\frac{\kappa+1}{2},\quad \Re(s_1+s_2)>\frac{\kappa+3}{2}
\end{align}
(using (i) and \cite[Theorem 3]{Mat02}).
As a generalization of \eqref{2-2}, define
\begin{align}\label{5-5}
F_{\pm}(s_1,s_2;\alpha;\omega;{\mathfrak A})=\sum_{l=1}^{\infty}
A_{s_1+s_2-1}(l;\pm\alpha;{\mathfrak A})
\Psi\left(s_2,s_1+s_2;\pm 2\pi i\omega l\right),
\end{align}
where
\begin{align}\label{5-6}
A_c(l;\pm\alpha;{\mathfrak A})=\sum_{mn=l}e^{\pm 2\pi im\alpha}a(n)n^c.
\end{align}
Then we have

\begin{thm}\label{thm5}
Assume {\rm (i)} and {\rm (ii)}.  Then the functions
$L_2(s_1,s_2;\alpha;\omega;{\mathfrak A})$ and
$F_{\pm}(s_1,s_2;\alpha;\omega;{\mathfrak A})$ can be continued
meromorphically to the whole space $\mathbb{C}^2$, and the functional equation
\begin{align}\label{5-7}
\lefteqn{L_2(s_1,s_2;\alpha;\omega;{\mathfrak A})}\\
&=\frac{\Gamma(1-s_1)\Gamma(s_1+s_2-1)}{\Gamma(s_2)\omega^{s_1+s_2-1}}
L(s_1+s_2-1,{\mathfrak A})\notag\\
&+\frac{\Gamma(1-s_1)}{\omega^{s_1+s_2-1}}\left\{
F_+(1-s_2,1-s_1;\alpha;\omega;{\mathfrak A})
+F_-(1-s_2,1-s_1;\alpha;\omega;{\mathfrak A})\right\}\notag
\end{align}
holds.
\end{thm}

In particular, applying Theorem \ref{thm5} to the special case $a(n)=1$ for only
one fixed $n$, and $a(n)=0$ for all other $n$, we obtain a functional equation for
the single series
\begin{align}\label{5-8}
\sum_{m=0}^{\infty}\frac{1}                                         
{(\alpha+m)^{s_1}(\alpha+m+n\omega)^{s_2}},
\end{align}
which gives the ``refinement'' (or the ``decomposition'') of Theorem \ref{thm5} 
in the sense of \cite[Remark 1.4]{CM}.   Note that \eqref{5-8} is a generalization
of the series
$$
\sum_{m=0}^{\infty}\frac{1}{(\alpha+m)^{s_1}(\beta+m)^{s_2}}
$$
($\beta\geq\alpha>0$) which was used in \cite{Mat03} and \cite{Mas04}.
Also Ram Murty and Sinha \cite{MurSin} studied analytic properties of this type of 
function and its generalizations.
 
Now we consider the case when $a(n)$'s are Fourier coefficients of modular forms.
Let $\mathcal{H}$ be the complex upper half plane, and
\begin{align}\label{5-9}
f(\tau)=\sum_{n=1}^{\infty}a(n)e^{2\pi i\tau n},
\end{align}
where $\tau\in\mathcal{H}$.    When $f(\tau)$ is a holomorphic cusp form of even
weight $\kappa$ with respect to the Hecke congruence subgroup $\Gamma_0(N)$,
it is well known that assumptions (i) and (ii) are satisfied.
Therefore, obviously, Theorem \ref{thm5} can be applied to this case.
In this case, however, we can show a different type of functional equation, 
using modular relations.   Let
\begin{align}\label{5-10}
\widetilde{f}(\tau)=(\sqrt{N}\tau)^{-\kappa}f\left(-\frac{1}{N\tau}\right)
\end{align}
and denote its Fourier expansion by 
$\widetilde{f}(\tau)=\sum_{n\geq 1}\widetilde{a}(n)e^{2\pi i\tau n}$.
Define
\begin{align}\label{5-11}
\lefteqn{
H_{2,N}^{\pm}(s_1,s_2;\alpha;\omega;\widetilde{f})}\\
&=
\sum_{m=1}^{\infty}\sum_{n=1}^{\infty}\widetilde{a}(n)e^{\mp 2\pi im\alpha}
m^{-s_1-s_2}\Psi\left(s_1+s_2,s_2;\pm\frac{2\pi in}{N\omega m}\right).
\notag
\end{align}
Write $L(s,f)$ and $L_2(s_1,s_2;\alpha;\omega;f)$  instead of
$L(s,{\mathfrak A})$ and $L_2(s_1,s_2;\alpha;\omega;{\mathfrak A})$,
respectively.    Then we have

\begin{thm}\label{thm6}
Assume {\rm (i)} and {\rm (ii)}, and $f(\tau)$ is as above.   Then the functions
$H_{2,N}^{\pm}(s_1,s_2;\alpha;\omega;\widetilde{f})$ can be continued
meromorphically to the whole space $\mathbb{C}^2$, and the functional equation
\begin{align}\label{5-12}
\frac{\Gamma(s_2)}{\Gamma(1-s_1)}
\lefteqn{L_2(s_1,s_2;\alpha;\omega;f)}\\                        
&\quad=\frac{ \Gamma(s_1+s_2-1)}{ \omega^{s_1+s_2-1}}        
L(s_1+s_2-1,f)\notag\\
&\quad+(2\pi)^{s_1+s_2-1}N^{-\kappa/2}\omega^{-\kappa}
 { \Gamma(\kappa-s_1-s_2+1)} \notag\\
&\quad\times\left\{e^{\pi i(1-s_1-s_2)/2}
H_{2,N}^+(-s_1,\kappa-s_2+1;\alpha;\omega;\widetilde{f})\right.\notag\\
&\quad\qquad\left.+e^{\pi i(s_1+s_2-1)/2}
H_{2,N}^-(-s_1,\kappa-s_2+1;\alpha;\omega;\widetilde{f})\right\}\notag 
\end{align}
holds.
\end{thm} 

\begin{rem}\label{rem3}
There is only one factor $a(n)$ on the numerator of \eqref{5-1}, \eqref{5-2}.
This is an unsatisfactory point; it is desirable to study the double series
with two factors (such as $a(m)b(n)$) on the numerator.   However our present
method cannot be applied to such a situation.
\end{rem}

In the remaining sections we give the proofs of Theorems \ref{thm5} and \ref{thm6}.
Since the proofs are straightforward generalization of the proofs described in
\cite{CM}, we just outline the argument briefly.
Denote $s_j=\sigma_j+it_j$ and $\theta=\arg\omega$.

%%%%%%%%%%%%%%%%%%%%%%%%%%%%%%%%%%%%%%%%%%%%%%%%%%%%%%%%%%%%%%%%%%%%%%%%%%%%%%%
\section{Sketch of the Proof of Theorem \ref{thm5}}
%%%%%%%%%%%%%%%%%%%%%%%%%%%%%%%%%%%%%%%%%%%%%%%%%%%%%%%%%%%%%%%%%%%%%%%%%%%%%%%

\textit{Step 1}.   The integral
\begin{align}\label{6-1}
\lefteqn{\Lambda(s_1,s_2;\alpha;\omega;{\mathfrak A})}\\
&=\int_0^{\infty}f(i\omega y)\int_0^{\infty}\frac{e^{2\pi(1-\alpha)(x+y)}}
{e^{2\pi(x+y)}-1}x^{s_1-1}y^{s_2-1}dxdy\notag
\end{align}
converges in the region
$\sigma_1>0$, $\sigma_2>(\kappa+1)/2$, $\sigma_1+\sigma_2>(\kappa+3)/2$, and
\begin{align}\label{6-2}
L_2(s_1,s_2;\alpha;\omega;{\mathfrak A})=
\frac{(2\pi)^{s_1+s_2}}{\Gamma(s_1)\Gamma(s_2)}
\Lambda(s_1,s_2;\alpha;\omega;{\mathfrak A}).
\end{align}
This can be shown by using \eqref{5-9} and the Taylor expansion
$$
\frac{e^{2\pi(1-\alpha)(x+y)}}{e^{2\pi(x+y)}-1}=
\sum_{m=0}^{\infty}e^{-2\pi m(x+y)-2\pi\alpha(x+y)}.
$$

\textit{Step 2}.  Let
\begin{align}\label{6-3}
h(z,\alpha)=\frac{e^{2\pi(1-\alpha) z}}{e^{2\pi z}-1}
-\frac{1}{2\pi z},
\end{align}
and divide the integral $\Lambda(s_1,s_2;\alpha;\omega;{\mathfrak A})$ as
\begin{align}\label{6-4}
\lefteqn{\Lambda(s_1,s_2;\alpha;\omega;{\mathfrak A})}\\
&=\int_0^{\infty}f(i\omega y)\int_0^{\infty}h(x+y,\alpha)x^{s_1-1}y^{s_2-1}dxdy
\notag\\
&\;+\int_0^{\infty}f(i\omega y)\int_0^{\infty}\frac{1}{2\pi(x+y)}
x^{s_1-1}y^{s_2-1}dxdy\notag\\ 
&=I_1+I_2,\notag
\end{align}
say.   The important point here is that the function $h(z,\alpha)$ is holomorphic at
$z=0$.    The idea of this type of decomposition goes back to Motohashi's short
note \cite{Mot85}, which influences \cite{KM91} (see Remark \ref{rem1}).

\textit{Step 3}.  By the beta integral formula, we find that
\begin{align}\label{6-5}
I_2=\frac{\Gamma(s_1)\Gamma(1-s_1)\Gamma(s_1+s_2-1)}{(2\pi)^{s_1+s_2}
\omega^{s_1+s_2-1}}L(s_1+s_2-1,{\mathfrak A})
\end{align}
in the region 
\begin{align}\label{6-6}
0<\sigma_1<1,\; \sigma_1+\sigma_2>(\kappa+3)/2.
\end{align}

\textit{Step 4}.   Using estimates of $h(z,\alpha)$ (\cite[p.1454]{KMT-IJNT}) and
$f(i\omega y)$ (\cite[Lemma 4.3.3]{Miy}), we find that the integral $I_1$ is
convergent in the region $\sigma_2>(\kappa+1)/2$.   Since \eqref{6-6} implies
$\sigma_2>(\kappa+1)/2$, now we know that the decomposition \eqref{6-4} is valid
in the region \eqref{6-6}.

\textit{Step 5}.   Changing the path of the inner integral of $I_1$ by the contour
$\mathcal{C}$ consisting of the half-line on the positive real axis from infinity
to a small positive number, a small circle counterclockwise round the origin,
and the other half-line on the positive real axis back to infinity.   Then
\begin{align}\label{6-7}
I_1=\frac{1}{e^{2\pi is_1}-1}I_3,
\end{align}
where
\begin{align}\label{6-8}                                                                
I_3=\int_0^{\infty}f(i\omega y)y^{s_2-1}\int_{\mathcal{C}}h(x+y,\alpha)               
x^{s_1-1}dxdy.                                                                          
\end{align}
Since $I_3$ is convergent in the region
\begin{align}\label{6-9}
\sigma_1<1,\; \sigma_2>(\kappa+1)/2,
\end{align}
and the right-hand side of \eqref{6-5} can be continued to the whole space
$\mathbb{C}^2$ (by assumption (ii)), now we see by \eqref{6-4} that
$\Lambda(s_1,s_2;\alpha;\omega;{\mathfrak A})$ is
(and hence $L_2(s_1,s_2;\alpha;\omega;{\mathfrak A})$ is) continued
meromorphically to the region \eqref{6-9}.

\textit{Step 6}.  Assume
\begin{align}\label{6-10}                                                               
\sigma_1<0,\; \sigma_2>(\kappa+1)/2                                                    
\end{align}
(a little smaller than \eqref{6-9}),
and change the path $\mathcal{C}$ on the right-hand side of \eqref{6-9} by
$$
\mathcal{C}_R=\{x=-y+(R+1/2)e^{i\varphi}\;|\;0\leq\varphi<2\pi\}
$$
($R\in\mathbb{N}$), and let $R\to\infty$.   Counting the residues, we have
\begin{align}\label{6-11}
I_3&=-\int_0^{\infty}f(i\omega y)y^{s_2-1}\sum_{m\neq 0}
e^{-2\pi im\alpha}\left(-y+im\right)^{s_1-1}dy\\
&=-i(I_{31}+I_{32}),\notag
\end{align}
where
\begin{align}\label{6-12}
I_{31}&=e^{\pi i(s_1-s_2-1)/2}\sum_{m=1}^{\infty}e^{-2\pi im\alpha}
m^{s_1+s_2-1}\\
&\times\int_0^{e^{i\pi/2}\infty} f\left(\omega mz\right)
z^{s_2-1}(z+1)^{s_1-1}dz,\notag
\end{align}
\begin{align}\label{6-13}                                                               
I_{32}&=e^{\pi i(3s_1+s_2-3)/2}\sum_{m=1}^{\infty}e^{2\pi im\alpha}   
m^{s_1+s_2-1}\\                                           
&\times\int_0^{e^{-i\pi/2}\infty} f\left(-\omega mz\right)   
z^{s_2-1}(z+1)^{s_1-1}dz.\notag                                                         
\end{align}

\textit{Step 7}.  Applying \eqref{5-9} to the right-hand side of \eqref{6-12} and
putting $mn=l$, we obtain
\begin{align}\label{6-14}
I_{31}&=e^{\pi i(s_1-s_2-1)/2}\Gamma(s_2)\\
&\times\sum_{l=1}^{\infty}A_{s_1+s_2-1}^0(l;-\alpha;{\mathfrak A})
\Psi\left(s_2,s_1+s_2;-2\pi i\omega l\right),\notag
\end{align}
where
\begin{align}\label{6-15}
A_c^0(l;\alpha;{\mathfrak A})=\sum_{mn=l}e^{2\pi im\alpha}m^c a(n).
\end{align}
Similarly we have
\begin{align}\label{6-16}                                                               
I_{32}&=e^{\pi i(3s_1+s_2-3)/2}\Gamma(s_2)\\                          
&\times\sum_{l=1}^{\infty}A_{s_1+s_2-1}^0(l;\alpha;{\mathfrak A})                      
\Psi\left(s_2,s_1+s_2;2\pi i\omega l\right).\notag                   
\end{align}
The interchanges of summation and integration in Steps 6 and 7 can be verified
similarly to the argument in \cite[Section 4]{CM}, in the region \eqref{6-10}.
Substituting \eqref{6-14} and \eqref{6-16} into \eqref{6-11} (and then \eqref{6-7}),
we obtain
\begin{align}\label{6-17}
\lefteqn{\frac{(2\pi)^{s_1+s_2}}{\Gamma(s_1)\Gamma(s_2)}I_1=
(2\pi)^{s_1+s_2-1}\Gamma(1-s_1)}\\
&\times\left\{e^{\pi i(1-s_1-s_2)/2}F_-^0(s_1,s_2;\alpha;\omega;
{\mathfrak A})+e^{\pi i(s_1+s_2-1)/2}F_+^0(s_1,s_2;\alpha;\omega;
{\mathfrak A})\right\},\notag
\end{align}
where
\begin{align}\label{6-18}
F_{\pm}^0(s_1,s_2;\alpha;\omega;{\mathfrak A})=\sum_{l=1}^{\infty}
A_{s_1+s_2-1}^0(l;\pm\alpha;{\mathfrak A})\Psi\left(s_2,s_1+s_2;
\pm 2\pi i\omega l\right).
\end{align}

\textit{Step 8}.   Using \eqref{2-6} and noting
\begin{align}\label{6-19}
A_c^0(l;\pm\alpha;{\mathfrak A})l^{-c}=A^{-c}(l;\pm\alpha;{\mathfrak A}),
\end{align}
we find that
\begin{align}\label{6-20}
F_{\pm}^0(s_1,s_2;\alpha;\omega;{\mathfrak A})=
\left(\pm 2\pi i\omega\right)^{1-s_1-s_2}
F_{\pm}(1-s_2,1-s_1;\alpha;\omega;{\mathfrak A}).
\end{align}
This is a generalization of \eqref{2-5}.
Substituting this into \eqref{6-17}, and combining with \eqref{6-5}, we now
arrive at formula \eqref{5-7}, in the region \eqref{6-10}.

\textit{Step 9}.   The final task is to show the meromorphic continuation of
the functions $L_2(s_1,s_2;\alpha;\omega;{\mathfrak A})$ and
$F_{\pm}(s_1,s_2;\alpha;\omega;{\mathfrak A})$.
The continuation of $F_{\pm}$ can be done by using the asymptotic expansion
\eqref{2-4} of the confluent hypergeometric function.   Then \eqref{5-7} gives
the continuation of $L_2$ to the whole space $\mathbb{C}^2$.   This completes the
proof of Theorem \ref{thm5}.

%%%%%%%%%%%%%%%%%%%%%%%%%%%%%%%%%%%%%%%%%%%%%%%%%%%%%%%%%%%%%%%%%%%%%%%%%%%%%%%
\section{Sketch of the proof of Theorem \ref{thm6}}
%%%%%%%%%%%%%%%%%%%%%%%%%%%%%%%%%%%%%%%%%%%%%%%%%%%%%%%%%%%%%%%%%%%%%%%%%%%%%%%%

Let $f(\tau)$ be as in the statement of Theorem \ref{thm6}.   Then 
$f(i\omega y)$ is of exponential decay as $y\to 0$, and hence the integral
\eqref{6-8} is convergent in the region $\sigma_1<1$ (wider than \eqref{6-9}).
Therefore, the formulas stated in Steps 6 and 7 in the preceding section are
now valid under the simple assumption $\sigma_1<0$.

Applying the modular relation \eqref{5-10} to \eqref{6-12} and \eqref{6-13}, and
changing the order of integration and summation, we obtain
\begin{align}\label{7-1}
I_{31}&=N^{-\kappa/2}\omega^{-\kappa}
e^{\pi i(s_1-s_2-1)/2}\\
&\times\Gamma(\kappa-s_1-s_2+1)
H_{2,N}^+(-s_1,\kappa-s_2+1;\alpha;\omega;\widetilde{f}),\notag
\end{align}
\begin{align}\label{7-2}                                                                
I_{32}&=N^{-\kappa/2}\omega^{-\kappa}                       
e^{\pi i(3s_1+s_2-3)/2}\\                                                                
&\times\Gamma(\kappa-s_1-s_2+1)
H_{2,N}^-(-s_1,\kappa-s_2+1;\alpha;\omega;\widetilde{f}).\notag        
\end{align}
Similarly to \cite[Section 6]{CM}, we can verify the above argument in the region
\begin{align}\label{7-3}
\sigma_1<0, \;\sigma_2<(\kappa-1)/2.
\end{align}
From \eqref{7-1}, \eqref{7-2} and \eqref{6-5}, we obtain \eqref{5-12} in the
region \eqref{7-3}.   Therefore the only remaining task is to prove the
meromorphic continuation of $H_{2,N}^{\pm}$.   We only discuss the case of
$H_{2,N}^+$.
   
Substituting the expression \eqref{2-1a} of the confluent hypergeometric function
into the right-hand side of \eqref{5-11}, putting $y=-i\eta$ and changing the order
of summation and integration, we obtain
\begin{align}\label{7-4}
H_{2,N}^+(s_1,s_2;\alpha;\omega;\widetilde{f})
&=\frac{-i}{\Gamma(s_1+s_2)}\int_0^{e^{i(\varphi+\pi/2)}\infty}
\sum_{m=1}^{\infty}e^{-2\pi im\alpha}m^{-s_1-s_2}\\
&\qquad\times\widetilde{f}\left(\frac{i\eta}{N\omega m}\right)
(-i\eta)^{s_1+s_2-1}(1-i\eta)^{-s_1-1}d\eta,\notag
\end{align}
where $\varphi$ satisfies the conditions $-\pi<\varphi<\pi$ and
$|\varphi+(\pi/2\theta)|<\pi/2$.
Putting
\begin{align}\label{7-5}                                                                
\widetilde{\mathcal{F}}(\tau,s,\alpha)=\sum_{m=1}^{\infty}e^{-2\pi im\alpha}            
m^{-s}\widetilde{f}\left(\frac{\tau}{Nm}\right),                                        
\end{align}
we obtain
\begin{align}\label{7-6}
\lefteqn{H_{2,N}^+(s_1,s_2;\alpha;\omega;\widetilde{f})}\\
&=\frac{-i}{\Gamma(s_1+s_2)}\int_0^{e^{i(\varphi+\pi/2)}\infty}
\widetilde{\mathcal{F}}\left(\frac{i\eta}{\omega},s_1+s_2,\alpha\right)
(-i\eta)^{s_1+s_2-1}(1-i\eta)^{-s_1-1}d\eta\notag\\
&=\frac{e^{i\varphi}}{\Gamma(s_1+s_2)}(-e^{i\varphi})^{s_1+s_2-1}\notag\\
&\quad\times\int_0^{\infty}\widetilde{\mathcal{F}}\left(-\frac{1}{\omega}
e^{i\varphi}\xi,s_1+s_2,\alpha\right)\xi^{s_1+s_2-1}(1+e^{i\varphi}\xi)^{-s_1-1}d\xi,
\notag
\end{align}
where $\xi=e^{-i(\varphi+\pi/2)}\eta$. 

As a generalization of \cite[Lemma 7.1]{CM}, we can show
\begin{align}\label{7-7}
\lefteqn{\int_0^{\infty}\widetilde{\mathcal{F}}\left(-\frac{1}{\omega}        
e^{i\varphi}\xi,s,\alpha\right)\xi^{u-1}d\xi}\\
&=\Gamma(u)\left(\frac{N\omega}{2\pi i}e^{-i\varphi}\right)^u
\phi(s-u,-\alpha)L(u,\widetilde{f})\notag
\end{align}
for $u\in\mathbb{C}$ with $(\kappa+1)/2<\Re u<\sigma-1$.   Therefore by the Mellin
inversion formula we have
\begin{align}\label{7-8}
\lefteqn{\widetilde{\mathcal{F}}\left(-\frac{1}{\omega}           
e^{i\varphi}\xi,s,\alpha\right)}\\
&=\frac{1}{2\pi i}\int_{(c)}\xi^{-u}\Gamma(u)\left(\frac{N\omega}{2\pi i}
e^{-i\varphi}\right)^u\phi(s-u,-\alpha)L(u,\widetilde{f})du,\notag
\end{align}
where $(\kappa+1)/2<c<\sigma-1$ and the path of integration is the vertical line
$\Re u=c$.
Substituting \eqref{7-8} into \eqref{7-6} and changing the order of integration
we obtain
\begin{align}\label{7-9}
H_{2,N}^+(s_1,s_2;\alpha;\omega;\widetilde{f})
&=\frac{e^{i\varphi}}{2\pi i\Gamma(s_1+s_2)}(-e^{i\varphi})^{s_1+s_2-1}
\int_{(c)}\Gamma(u)\left(\frac{N\omega_2}{2\pi i\omega_1}       
e^{-i\varphi}\right)^u\\
&\times\phi(s-u,-\alpha)L(u,\widetilde{f})\int_0^{\infty}\xi^{s_1+s_2-1-u}
(1+e^{i\varphi}\xi)^{-s_1-1}d\xi du.\notag
\end{align}
Using the beta integral formula we can show that the inner integral on the
right-hand side of \eqref{7-9} is
$$
=(e^{-i\varphi})^{s_1+s_2-u}\frac{\Gamma(u-s_2+1)\Gamma(s_1+s_2-u)}
{\Gamma(s_1+1)},
$$
and hence
\begin{align}\label{7-10}
H_{2,N}^+(s_1,s_2;\alpha;\omega;\widetilde{f})
&=\frac{e^{i\varphi}(-e^{i\varphi})^{s_1+s_2-1}}{2\pi i\Gamma(s_1+s_2)
\Gamma(s_1+1)}\int_{(c)}\Gamma(u)\Gamma(u-s_2+1)\Gamma(s_1+s_2-u)\\
&\times\left(\frac{N\omega}{2\pi i}                                 
e^{-i\varphi}\right)^u\phi(s-u,-\alpha)L(u,\widetilde{f})du.\notag
\end{align}
Finally, modifying the path of integration on the right-hand side of \eqref{7-10}
in the same way as in \cite[Section 7]{CM}, we can show the meromorphic 
continuation of $H_{2,N}^+$.   The case of $H_{2,N}^-$ is similar.
The proof of Theorem \ref{thm6} is complete.

%%%%%%%%%%%%%%%%%%%%%%%%%%%%%%%%%%%%%%%%%%%%%%%%%%

\end{document}